\documentclass[a4paper,12pt]{article}

\usepackage{graphicx}
\usepackage[centertags]{amsmath}
\usepackage{graphicx}
\usepackage{a4wide}
\usepackage{amsfonts}
\usepackage{amsbsy}
\usepackage{amssymb}
\usepackage{amscd}

\def\prod#1{\langle{#1}\rangle}

\def\squarebox#1{\hbox to #1{\hfill\vbox to #1{\vfill}}}
\newcommand{\qed}{\hspace*{\fill}
\vbox{\hrule\hbox{\vrule\squarebox{.667em}\vrule}\hrule}\smallskip}
\newtheorem{teorema}{Theorem}[section]
\newtheorem{lema}[teorema]{Lemma}
\newtheorem{corolario}[teorema]{Corollary}

\newenvironment{prova}{\noindent {\bf Proof:}}{\hfill $\qed $ \newline}

\newcommand{\n}{{\mathbb N}}
\newcommand{\Z}{{\mathbb Z}}

\newcommand{\g}{{\mathfrak g}}

\renewcommand{\a}{{\mathfrak a}}
\renewcommand{\n}{{\mathfrak n}}
\newcommand{\z}{{\mathfrak z}}

\newcommand{\p}{{\mathfrak p}}

\newcommand{\Ad}{{\rm Ad}}

\newcommand{\T}{\Theta}

\newcommand{\D}{\Delta}
\newcommand{\F}{{\mathbb F}}

\newcommand{\cl}{{\rm cl}}

\begin{document}

\title{A note on the Bruhat decomposition \\of semisimple Lie groups}

\author{Lucas Seco\thanks{Supported by FAPESP grant
n$^{\protect\underline{\circ }}$ 07/52390-6}}

\maketitle

\begin{abstract}
Let a split element of a connected semisimple Lie group act on one of its flag
manifolds. We prove that each connected set of fixed points of this
action is itself a flag manifold. With this we can obtain the 
generalized Bruhat decomposition of a semisimple Lie
group by entirely dynamical arguments.
\end{abstract}



\setcounter{section}{-1}

\section{Introduction}

Let a split element $h$ of a connected semisimple Lie group $G$ act on one of its flag
manifolds $\F_\T$ (notation of semisimple Lie groups is recalled in Section
\ref{preliminar-lie}). We prove that each connected set of fixed points of this
action is itself a flag manifold, but a flag manifold of a semisimple Lie subgroup of $G$. This generalizes directly the fact that each connected fixed point set of a diagonable matrix acting on a projective space is given by a projective subspace.  Apart for being interesting in itself, this result also allows us to obtain generalized Bruhat decomposition of a semisimple Lie group by entirely dynamical arguments, as we explain below.

Standard textbooks on semisimple Lie groups \cite{Helgason, w} prove
by entirely algebraic arguments what we will call the regular Bruhat
decomposition of a connected semisimple Lie group $G$, namely
$$
G = \coprod_{w \in W} P w P = \coprod_{w \in W} N^+ w P,
$$
where $P$ is the minimal parabolic and $W$ the Weyl group of $G$.
This decomposition corresponds to the regular
Bruhat decomposition of the maximal flag manifold $\F = \Ad(G)\p$ of
$G$, given by
$$
\F = \coprod_{w \in W} P w \p = \coprod_{w \in W} N^+ w \p,
$$
which can be seen as the decomposition of $\F$
into unstable manifolds of the action of a split-regular element $h
\in A^+$: see Section 3 of \cite{dkv} for a proof of this by
entirely dynamical arguments (this is recalled at the end of Section
\ref{preliminar-lie}).
From this regular Bruhat decomposition on the maximal flag manifold
$\F$ one readily obtains the regular Bruhat decomposition on the
partial flag manifolds $\F_\T = \Ad(G)\p_\T$ given by
$$
\F_\T = \coprod_{w \in W/W_\T} P w \p_\T = \coprod_{w \in W /W_\T}
N^+ w \p_\T,
$$
the argument goes as follows.  Projecting the regular Bruhat
decomposition of $\F$ onto $\F_\T$ one needs only to show the
disjointedness of the above decomposition. If the unstable manifolds
$N^+ w \p_\T$ and $N^+ s \p_\T$ meet, for $s, w \in W$, then there
exists $n \in N^+$ such that $w\p_\T = n s\p_\T$. Taking the regular
element $h \in A^+$ we have for $k \in \Z$ that $w\p_\T$ is a fixed
point so that
$$
w\p_\T = h^{-k} w\p_\T = h^{-k} n s\p_\T \to s\p_\T,
$$
when $k \to +\infty$.  It follows that $w\p_\T = s\p_\T$, so that
$s^{-1} w\p_\T = \p_\T$ which, by the Iwasawa decomposition of
$P_\T$, implies that $s^{-1} w \in K_\T \cap M^*$ so that $s^{-1} w
\in W_\T$, that is, $w \in sW_\T$, as claimed.  The corresponding
decomposition in $G$ is the regular Bruhat decomposition
$$
G = \coprod_{w \in W/W_\T} P w P_\T = \coprod_{w \in W /W_\T} N^+ w
P_\T.
$$

Usually much harder to obtain is what we will call the generalized Bruhat decomposition of $G$, given by
$$
G = \coprod_{w \in W_\D \backslash W/W_\T} P_\D w P_\T.
$$
To show this \cite{w} appeals to the use of Tits Systems (see Section
1.2 of \cite{w}). Dynamically, this corresponds to the
decomposition of the flag manifold $\F_\T$ into unstable manifolds
of the action of an irregular split-element $h \in \cl A^+$, $h =
\exp(H)$,
$$
\F_\T = \coprod_{w \in W_H \backslash W/W_\T} P_H w \p_\T =
\coprod_{w \in W_H \backslash W/W_\T} N^+_H Z_H w \p_\T,
$$
where $\D$ is the set of simple roots which annihilate $H$.
In this case the fixed points of $h$ degenerate into fixed point
manifolds $Z_H w \p_\T$.  To show that these unstable manifolds are disjoint we can argue as above to get rid of the unstable part $N^+_H$ so that the only difficulty
is to show that fixed point manifolds are disjoint when we take $ w \in
W_H \backslash W/W_\T$.  At this point of the argument \cite{dkv} appeals to a general theorem of Borel-Tits (see Proposition 1.3 of \cite{dkv}).

In this note we show the disjointedness of the above fixed point manifolds as a byproduct of 
showing that each of these fixed point manifold is itself
equivariantly diffeomorphic to a flag manifold.  With this, one can obtain the
generalized Bruhat decomposition of a semisimple Lie
group by entirely dynamical arguments: one follows Section 3 of \cite{dkv} to prove the regular Bruhat decomposition and then uses the result of this article to prove the generalized Bruhat decomposition.

In the first section we recall notation and preliminar results on
semisimple Lie theory.  In the second section we prove the main
result Theorem \ref{teo:fix-flag} and deduce from it that the fixed
point manifolds are disjoint.

\section{Preliminaries on Semi-simple Lie Theory\label{preliminar-lie}}
For the theory of semi-simple Lie groups and their flag manifolds we
refer to Duistermat-Kolk-Varadarajan \cite{dkv}, Helgason
\cite{Helgason} and Warner \cite{w}. To set notation let $G$ be a
connected noncompact semi-simple Lie group  with  Lie algebra
$\frak{g}$. We assume throughout that $G$ has finite center. Fix a
Cartan involution $\theta $ of $\frak{g}$ with Cartan decomposition
$\frak{g}=\frak{k}\oplus \frak{s}$. The form $B_{\theta }\left(
X,Y\right) =-\langle X,\theta Y\rangle $, where $\langle \cdot
,\cdot \rangle $ is the Cartan-Killing form of $\frak{g}$, is an
inner product.

Fix a maximal abelian subspace $\frak{a}\subset \frak{s}$ and a Weyl
chamber $\frak{a}^{+}\subset \frak{a}$. We  let $\Pi $ be the set of
roots of $\frak{a}$,  $\Pi ^{+}$ the positive roots corresponding to
$\frak{a}^{+}$, $\Sigma $ the set of simple roots in $\Pi ^{+}$ and
$\Pi^- = - \Pi^+$ the negative roots.  For a root $\alpha \in \Pi$
we denote by $H_\alpha \in \frak{a}$ its coroot so that
$B_\theta(H_\alpha,H) = \alpha(H)$ for all $H \in \frak{a}$. The
Iwasawa decomposition of the Lie algebra $\frak{g}$ reads
$\frak{g}=\frak{k}\oplus \frak{a}\oplus \frak{n}^{\pm}$ with
$\frak{n}^{\pm}=\sum_{\alpha \in \Pi ^{\pm}}\frak{g}_{\alpha }$
where $\frak{g}_{\alpha }$ is the root space associated to $\alpha
$.  As to the global decompositions of the group we write  $G=KS$
and $G=KAN^\pm$ with $K=\exp \frak{k}$, $S=\exp \frak{s}$, $A=\exp
\frak{a}$ and $N^{\pm}=\exp \frak{n}^{\pm}$.

The Weyl group $W$ associated to $\frak{a}$ is the finite group
generated by the reflections over the root hyperplanes $\alpha = 0$
in $\frak{a}$, $\alpha \in \Pi$. $W$ acts on $\frak{a}$ by
isometries and can be alternatively be given as $W=M^{*}/M$ where
$M^{*}$ and $M$ are the normalizer and the centralizer of $A$ in
$K$, respectively. We write $\frak{m}$ for the  Lie algebra of $M$.
There is an unique element $w^- \in W$ which takes the simple roots
$\Sigma $ to $-\Sigma$, $w^-$ is called the principal involution of
$W$.

Associated to  a subset of simple roots $\Theta \subset \Sigma $
there are several Lie  algebras and groups (cf.  \cite{w}, Section
1.2.4): We write $\frak{g}\left( \Theta \right) $ for the
(semi-simple)  Lie subalgebra generated by $\frak{g}_{\alpha }$,
$\alpha \in \Theta $, and put $\frak{k}(\Theta) = \frak{g}(\Theta)
\cap \frak{k}$, $\frak{a} \left( \Theta \right) =\frak{g}\left(
\Theta \right) \cap \frak{a}$, and $\frak{n}^{\pm }\left( \Theta
\right) =\frak{g}\left( \Theta \right) \cap \frak{n}^{\pm }$. The
simple roots of $\frak{g}(\Theta)$ are given by $\Theta$, more
precisely, by restricting the functionals of $\Theta$ to
$\frak{a}(\Theta)$. The coroots $H_\alpha$, $\alpha \in \Theta$,
form a basis for $\frak{a} \left( \Theta \right)$.  Let $G\left(
\Theta \right) $ and $K(\Theta)$ be the connected groups with Lie
algebra $\frak{g}\left( \Theta \right)$ and $\frak{k}\left( \Theta
\right)$, respectively. Then $G(\Theta)$ is a connected semisimple
Lie group with finite center.  Let $A\left( \Theta \right) =\exp
\frak{a}\left( \Theta \right) $, $N^{\pm }\left( \Theta \right)
=\exp \frak{n}^{\pm }\left( \Theta \right) $.  We have the Iwasawa
decomposition $G(\Theta) = K(\Theta) A(\Theta) N^\pm(\Theta)$. Let
$\frak{a}_{\Theta }=\{H\in \frak{a}:\alpha(H) =0,\, \alpha \in
\Theta \}$  be the orthocomplement of $\frak{a}(\Theta)$ in
$\frak{a}$ with respect to the $B_\theta$-inner product and put
$A_{\Theta }=\exp \frak{a}_{\Theta }$. The subset $\Theta $ singles
out the subgroup $W_\T$ of the Weyl group which acts trivially
$\frak{a}_\T$.  Alternatively $W_{\Theta }$ can be given as the
subgroup generated by the reflections with respect to the roots
$\alpha \in \Theta$.  The restriction of $w \in W_\T$ to
$\frak{a}(\T)$ furnishes an isomorphism between $W_\T$ and the Weyl
groyp $W(\T)$ of $G(\T)$

The standard parabolic subalgebra of type $\T \subset \Sigma$ with
respect to chamber $\frak{a}^+$ is defined by
\[
\frak{p}_{\Theta }=\frak{n}^{-}\left( \Theta \right) \oplus {\frak
m}\oplus {\frak a}\oplus {\frak n}^{+}.
\]
The corresponding standard parabolic subgroup $P_{\Theta }$ is the
normalizer of $\frak{p}_{\Theta }$ in $G$.  It has the Iwasawa
decomposition $P_{\Theta } = K_\T A N^+$. The empty set $\Theta
=\emptyset $ gives the minimal parabolic subalgebra $\frak{p} =
{\frak m}\oplus {\frak a}\oplus {\frak n}^{+}$ whose minimal
parabolic subgroup $P=P_\emptyset $ has Iwasawa decomposition
$P=MAN^+$.  Let $\D \subset \T \subset \Sigma$.  Then $P_\D \subset
P_\T$, also we denote by $P(\T)_\D$ the parabolic subgroup of
$G(\T)$ of type $\D$.

We let  $Z_{\Theta }$ be the  centralizer of $\frak{a}_{\Theta }$ in
$G$ and $K_{\Theta }=Z_{\Theta }\cap K$.  We have that $K_\T$
decomposes as $K_\T = M K(\T)$ and that $Z_\T$ decomposes as
$Z_{\Theta }=MG(\Theta )A_{\Theta }$ which implies that
$Z_{\Theta }=K_{\Theta }AN(\Theta )$ is the Iwasawa decomposition of
$Z_{\Theta }$ (which is a reductive Lie group).  Let $\Delta \subset
\Sigma$, then\footnote{Using that $\frak{a}_\Theta =
\frak{a}(\Theta)^\perp$ and that $\frak{a}(\Theta\cap\Delta) =
\frak{a}(\Theta)\cap\frak{a}(\Delta)$ this follows by taking perp on
both sides and using that $(V+W)^\perp = V^\perp \cap W^\perp$ for
$V,W$ linear subspaces. }
$\frak{a}_{\Theta\cap\Delta} = \frak{a}_{\Theta} +
\frak{a}_{\Delta}$.  Thus it follows that $Z_{\Theta \cap \Delta} =
Z_\Theta \cap Z_\Delta$, $K_{\Theta \cap \Delta} = K_\Theta \cap
K_\Delta$ and $P_{\Theta \cap \Delta} = P_\Theta \cap P_\Delta$.
For $H \in \a$ we denote by $Z_H$, $W_H$ etc.\ the centralizer of
$H$ in $G$, $W$ etc., respectively.  When $H\in \cl\frak{a}^+$ we
put
$$
\Theta(H)=\{\alpha \in \Sigma :\alpha (H)=0\},
$$
and we have that $Z_{H}=Z_{\Theta(H)}$, $K_{H}=K_{\Theta(H)}$,
$N^+_H = N^+(\T(H))$ and $W_{H}=W_{\Theta(H)}$.

Let $\frak{n}^{\pm}_\T= \sum_{\alpha \in \Pi^{\pm} - \langle \T
\rangle} \frak{g}_{\alpha }$ and $N^{\pm}_\T =
\exp(\frak{n}^{\pm}_\T)$.  Then $N^{\pm}$ decomposes as $N^{\pm} =
N(\T)^{\pm} N_\T^{\pm}$ where $N(\T)^{\pm}$ normalizes $N_\T^{\pm}$
and $N(\T)^{\pm} \cap N_\T^{\pm} = 1$. We have that $\mathfrak{g} =
\frak{n}^{-}_\T \oplus \frak{p}_\T$, that $N^-_\T \cap P_\T = 1$ and
also that $P_{\Theta }$ is the normalizer of $\frak{n}^{+}_\T$ in
$G$. $P_\T$ decomposes as $P_\T = Z_\T N^+_\T$, where $Z_\T$
normalizes $N^+_\T$ and $Z_\T \cap N^+_\T = 1$.
We write $\frak{p}_{\Theta }^{-}=\theta (\frak{p}_{\Theta })$ for
the parabolic subalgebra opposed to $\frak{p}_{\Theta }$. It is
conjugate to the parabolic subalgebra $\frak{p}_{\Theta ^{*}}$ where
$\Theta ^{*}=-(w^-)\Theta $ is the dual to $\Theta$ and $w^-$ is the
principal involution of $W$. More precisely, $\frak{p}_{\Theta
}^{-}=k\frak{p}_{\Theta^{*}}$ where $k \in M^{*}$ is a
representative of $w^-$. If $P_{\Theta }^{-}$ is the parabolic
subgroup associated to $\frak{p}_{\Theta }^{-}$ then $Z_{\Theta
}=P_{\Theta }\cap P_{\Theta }^{-}$ and $P_{\Theta }^{-}=Z_\T
N^-_\T$, where $Z_\T$ normalizes $N^-_\T$ and $Z_\T \cap N^-_\T =
1$.

The flag manifold  of type $\Theta $ is the orbit $\Bbb{F}_{\Theta
}={\rm Ad}(G)\frak{p}_\T$, which identifies with the homogeneous
space $G/P_{\Theta}$.  Since the center of $G$ normalizes
$\frak{p}_\T$, the flag manifold depends only on the Lie algebra
$\frak{g}$ of $G$.  The empty set $\Theta =\emptyset $ gives the
maximal flag manifold $\Bbb{F}=\Bbb{F}_\emptyset $.  If $\D \subset
\T$ then there is a $G$-equivariant projection $\Bbb{F}_\D \to
\Bbb{F}_\T$ given by $g\frak{p}_\D \mapsto g\frak{p}_\T$, $g \in G$.

Some subalgebras of $\g$ which are defined by the choice of a Weyl
chamber of $\a$ and a subset of the associated simple roots can be
defined alternatively by the choice of an element $H \in \a$ as
follows. First note that the eigenspaces of ${\rm ad}(H)$ in $\g$
are the weight spaces $\g_\alpha$, and that the centralizer of $H$
in $\g$ is given by $ \z_H = \sum\{ \g_\alpha:\, \alpha(H) = 0 \}$,
where the sum is taken over $\alpha \in \a^*$. Now define the
negative and positive nilpotent subalgebras of type $H$ given by
$$
\n^-_H = \sum\{ \g_\alpha:\, \alpha(H) < 0 \}, \qquad \n^+_H =
\sum\{ \g_\alpha:\, \alpha(H) > 0 \},
$$
and the parabolic subalgebra of type $H$ which is given by
$$
 \p_H = \sum\{
\g_\alpha:\, \alpha(H) \geq 0 \},
$$
where in all cases $\alpha$ runs through all the roots $\Pi$. Then
we have that
$$
\g = \n^-_H \oplus \z_H \oplus \n^+_H\qquad\mbox{and}\qquad \p_H =
\z_H\oplus \n^+_H,
$$
and that
$$
\n^\pm_{wH} = w\n^\pm_H \qquad \p_{wH} = w\p_H.
$$
Define the flag manifold of type $H$ given by the orbit
$$
\F_H = \Ad(G)\p_H.
$$
Now choose a chamber $\a^+$ of $\a$ which contains $H$ in its
closure, consider the simple roots $\Sigma$ associated to $\a^+$ and
take $\T(H) \subset \Sigma$.  Since a root $\alpha \in \T(H)$ if,
and only if, $\alpha|_{\a_{\T(H)}} = 0$, we have that
$$
\z_H = \z_{\T(H)},\quad \n^\pm(H) = \n^\pm_{\T(H)}, \quad \p_H =
\p_{\T(H)}.
$$
So it follows that
$$
\F_H = \F_{\T(H)},
$$
and that the isotropy of $G$ in $\p_H$ is $P_{\T(H)} = K_{\T(H)} A
N^+ = K_H A N^+$, since $K_{\T(H)} = K_H$.  In particular we have
that
$$
w\p_{\T(H)} = w \p_H = \p_{wH}.
$$
We note that we can proceed reciprocally. That is, if $\a^+$ and
$\T$ are given, we can choose an $H \in \cl\a^+$ such that $\T(H) =
\T$ and describe the objects that depend on $\a^+$ and $\T$ by $H$
(clearly, such an $H$ is not unique.)

We recall the dynamics of a split element of $G$ acting in the flag
manifold $\F_\T$ (see Section 3 of \cite{dkv}).  Let the split
element $H \in \cl\frak{a}^+$ act in $\F_\T$ by
\begin{equation}
t \cdot  b = \exp(t H)  b,\quad  b \in \F_\T.
\end{equation}
The connected sets of fixed point of the $H$-action are parametrized
by $w \in W$ and are given by
\begin{equation}\label{def:fixed-points-H}
{\rm fix}(H,w) =  Z_{H} w \frak{p}_\T  =  K_{H} w \frak{p}_\T ,
\end{equation}
so that they are in bijection with the double coset $W_{H}
\backslash W / W_\T$. It follows from $Z_H = Z_{\T(H)}$ and from the
decomposition $Z_{\Theta }=G(\Theta ) MA_{\Theta }$ that we have
$$
{\rm fix}(H,w) =  G(\T(H)) w \frak{p}_\T.
$$
Each $w$-fixed point set has stable/unstable manifold given
respectively by
\begin{equation}\label{def:stable-manifold-H}
{\rm st}(H,w) =  N^\pm_{\T(H)} {\rm fix}(H,w)  =  P_{\T(H)}^\pm w
\frak{p}_\T .
\end{equation}
This $H$-action decomposes $\F_\T$ in the disjoint union of
stable/unstable manifolds
\begin{equation}\label{def:bruhat-H}
\F_\T = \coprod_{W_{H} \backslash W / W_\T} P_{\T(H)}^\pm w
\frak{p}_\T,
\end{equation}
this is known as the Bruhat decomposition of $\F_\T$.  It follows
from these considerations that the dynamics of the $H$-action in
$\F_\T$ depends only on the set of roots ${\T(H)}$ which annihilate
$H$.

\section{Fixed points as flag manifolds}

Let $\pi_\T: \a \to \a(\T)$ be the orthogonal projection parallel to
$\a_\T$.

\begin{lema} The following assertions are true.
\begin{enumerate}
\item The projection by $\pi_\T$ of a regular element of $\a$ is a regular element of $\a(\T)$.
\item The projection by $\pi_\T$ of a chamber in $\a$ is contained inside a chamber of $\a(\T)$.
\item For $w \in W$ denote by $\a(\T)^w$ the chamber of $\a(\T)$ which contains
the projection by $\pi_\T$ of the chamber $w\a^+$.  Then the
nilpotent subalgebras $\n^+$ e $\n(\T)^w$ w.r.t.\ to the chambers
$\a^+$ and $\a(\T)^w$ satisfy
$$ \n(\T)^w \subset w \n^+. $$
\end{enumerate}
\end{lema}
\begin{prova}
We first observe that for $\alpha \in \T$, we have that
$\alpha|_{\a_\T} = 0$ so that for $H \in \a$ we have
$\alpha(\pi_\T(H)) = \alpha(H)$.  From this it follows that if $H$
is regular in $\a$, then $\pi_\T(H)$ is regular in $\a(\T)$, which
proves the first item. For the second item we observe that the
projection of a chamber of $\a$ is a convex set of $\a(\T)$ which,
by the first item, consists of regular elements of $\a(\T)$, and
hence it is contained in a chamber of $\a(\T)$. For the third item
let $\alpha \in \T$. If $\alpha > 0$ in $\a(\T)^w$ then $\alpha > 0$
in $\pi_\T(w\a^+)$ and hence, by the first remark of the proof, we
have that $\alpha > 0$ in $w\a^+$.  It follows that
$$
\n(\T)^w = \sum \{\g_\alpha:\, \alpha|_{\a(\T)^w} > 0, \alpha \in
\prod{\T} \} \subset \sum \{\g_\alpha:\, \alpha|_{w\a^+} > 0, \alpha
\in \Pi\} = w \n^+,
$$
as desired.
\end{prova}

\begin{teorema}\label{teo:fix-flag}
Let $X \in \cl\a^+$, $\T \subset \Sigma$, $w \in W$. Consider $\D =
\T(X)$ and $H_\T \in \cl\a^+$ such that $\T(H_\T) = \T$. Then the
map
$$
{\rm fix}(X,w)_\T \to \F_{\pi_\D(wH_\T)}(\g(\D)),\quad
g \p_{wH_\T} \mapsto g \p_{\pi_\D(wH_\T)}, \quad g \in G(\D)
$$
is a well defined $G(\D)$-equivariant diffeomorphism.
\end{teorema}
\begin{prova}
Since we have ${\rm fix}(X,w)_\T = G(\D) \p_{wH_\T}$ it follows that
the above map, let us call it $\psi$, is defined in all of its
domain and it is $G(\D)$-equivariant. It remains to prove that
$\psi$ is well defined and is injective both of which will follow if
we show that the isotropy of $\p_{wH_\T}$ in $G(\D)$ coincides with
the isotropy of $\p_{\pi_\D(wH_\T)}$ in $G(\D)$.  For this, let
$\a(\D)^w$ the chamber of $\a(\D)$ which contains the projection by
$\pi_\D$ of the chamber $w\a^+$. Consider the Iwasawa decomposition
of $G$ and $P_{w H_\T}$ w.r.t.\ to the chamber $w\a^+$, and the
Iwasawa decomposition of $G(\D)$ w.r.t.\ to the chamber $\a(\D)^w$
$$
G = KAwN^+w^{-1}, \quad P_{w H_\T}= K_{w H_\T}AwN^+w^{-1}, \quad
G(\D) = K(\D)A(\D)N(\D)^w,
$$
where, by item (3) of the previous Lemma, we have that
$$ N(\D)^w \subset wN^+w^{-1}.$$
Thus, by the uniqueness of the Iwasawa decomposition of $G$, it
follows that the isotropy of $\p_{wH_\T}$ in $G(\D)$ is given by
$$
G(\D) \cap P_{w H_\T} = (K(\D) \cap K_{w H_\T}) A(\D) N(\D)^w.
$$
The first term in the right hand side can be written as
$$
K(\D) \cap K_{w H_\T} = K(\D)_{w H_\T} = K(\D)_{\pi_\D(w H_\T)},
$$
where in the last equality we used that $K(\D)$ already centralizes
$\a_\D$.  Since $\pi_\D(w H_\T)$ lies in the closure of the chamber
$a(\T)^w$, it follows that
$$
G(\D) \cap w P_\T w^{-1} = K(\D)_{\pi_\D(w H_\T)} A(\D) N(\D)^w =
P_{\pi_\D(w H_\T)}(\D),
$$
which is precisely the isotropy of $\p_{\pi_\D(wH_\T)}$ in $G(\D)$.
It is then immediate that the inverse of $\psi$ is given by $g
\p_{\pi_\D(wH_\T)} \mapsto g \p_{wH_\T}$, $g \in G(\D)$, which shows
that $\psi$ is a diffeomorphism.
\end{prova}

\begin{corolario}
If ${\rm fix}(X,w')_\T \cap {\rm fix}(X,w)_\T \neq \emptyset$ then
$w' \in W_XwW_\T$.
\end{corolario}
\begin{prova}
Here we will adopt the notation of Theorem \ref{teo:fix-flag},
denoting by $\psi$ the diffeomorphism of that theorem. If ${\rm
fix}(X,w')_\T \cap {\rm fix}(X,w)_\T \neq \emptyset$ then there
exists $g \in G(\D)$ such that $w'\p_\T = g w \p_\T$. Take a regular
$h \in A(\D)$, using the $G(\D)$-equivariance of $\psi$ for $k \in
\Z$ we have that
$$
w' \p_\T = h^k w' \p_\T = h^k g w \p_\T = h^k g\p_{wH_\T}  =
\psi^{-1}( h^k g \p_{\pi_\D(wH_\T)}) = (*).
$$
By the regular Bruhat decomposition of the flag manifold
$\F_{\pi_\D(wH_\T)}(\g(\D))$ (cf.\ Theorem \ref{}), letting $k \to
\infty$ we have that there exists $s \in W(\D) = W_\D = W_X$ such
that
$$
(*) \to \psi^{-1}( s\p_{\pi_\D(wH_\T)}) = s\p_{wH_\T} = sw \p_\T.
$$
It follows that $w^{-1} s^{-1} w' \p_\T = \p_\T$, so that $w^{-1}
s^{-1} w' \in M^* \cap K_\T$, which implies that $w^{-1} s^{-1} w
\in W_\T$. Hence $w \in swW_\T \subset W_XwW_\T$, as desired.
\end{prova}


\end{document}